\documentclass[12pt]{article}
\usepackage{amsmath,amssymb}

\def\seq#1#2#3{#1_{#2},\,\ldots,#1_{#3}}
\newcommand{\Var}{\ensuremath{\mathcal{V}_{\mathbb{C}}}}

\def\vv{{\underline{v}}}
\def\tt{{\underline{t}}}
\def\aa{\underline{a}}
\def\ww{{\underline{w}}}
\def\mm{\underline{m}}

\def\1{\underline{1}}

\def\P{\Bbb P}

\def\n{{\underline{{\bf n}}}}
\def\LL{\Bbb L}
\def\Z{\Bbb Z}

\def\DD{\Bbb D}
\def\K{\Bbb K}
\def\C{\Bbb C}
\def\A{\Bbb A}

\def\OO{{\cal O}}
\def\X{{\cal X}}
\def\D{{\cal D}}
\def\Init{\mathsf{Init\,}}

\newtheorem{theorem}{Theorem}

\newtheorem{lemma}{Lemma}
\newtheorem{proposition}{Proposition}

\newenvironment{definition}
{\smallskip\noindent{\bf Definition\/}:}{\smallskip\par}

\newenvironment{example}
{\smallskip\noindent{\bf Example\/}.}{\smallskip\par}

\newenvironment{remark}
{\smallskip\noindent{\bf Remark\/}.}{\smallskip\par}

\newenvironment{proof}
{\noindent{\bf Proof\/}.}{{ $\Box$}\smallskip\par}

\title{Multi-index filtrations and motivic Poincar\'e series.
\footnote{Math. Subject Class. 14H20, 32S99}
}
\author{A.~Campillo
\and F.~Delgado \thanks{First two authors were partially
supported by the grant MCyT BFM2001-2251. Address:
University of Valladolid, Dept. of Algebra, Geometry and Topology,
47005 Valladolid, Spain.
E-mail: campillo\symbol{'100}cpd.uva.es, fdelgado\symbol{'100}agt.uva.es}
\and S.M.~Gusein-Zade \thanks{Partially supported by the grants
RFBR--04--01--00762,
NSh--1972.2003.1.
The author is thankful to the University of Valladolid for hospitality.
Address: Moscow State University,
Faculty of Mathematics and Mechanics, Moscow, GSP-2, 119992, Russia.
E-mail: sabir\symbol{'100}mccme.ru}}
\date{}
\begin{document}
\sloppy
\def\eps{\varepsilon}

\maketitle

\begin{abstract}
The Poincar\'e series of a multi-index filtration on the ring of germs of
functions can be written as a certain integral with respect to the Euler
characteristic over the projectivization of the ring. Here this integral
is considered
with respect to the generalized Euler characteristic with values in the
Grothendieck ring of varieties. For the filtration defined by orders of
functions on the components of a plane curve singularity $C$ and for the so
called divisorial filtration for a modification of $(\C^2,0)$ by a sequence of
blowing-ups there are given  formulae for this integral in terms of an
embedded
resolution of the germ $C$ or in terms of the modification respectively.
The generalized Euler characteristic of the extended semigroup
corresponding to the divisorial filtration is computed.
\end{abstract}

\section*{Introduction}

Let $(C,0)\subset (\C^n,0)$ be a germ of a reduced analytic curve, let
$C =\bigcup\limits^{r}_{k=1} C_k$ be its decomposition into irreducible
components, and let
$\varphi_k : (\C,0)\to (\C^n,0)$, $k=1,\ldots, r$, be parametrizations of the
components $C_k$ of the curve $C$; i.e. germs of analytic maps such that
$\mbox{Im}\, \varphi_k =C_k$ and $\varphi_k$ is an isomorphism between $\C$ and
$C_k$ outside of the origin (in a neighbourhood of it).
For a germ $g\in \OO_{\C^n,0}$ let $v_k=v_k(g)$ and $a_k=a_k(g)$ be the power of
the leading term  and the coefficient at it in the power series decomposition of
the germ $g\circ \varphi_k: (\C,0)\to \C$,
$g\circ \varphi_k (\tau) = a_k \tau^{v_k} + \mbox{terms of higher degree}$,
$a_k\neq 0$.  If $g\circ \varphi_k (\tau) \equiv 0$, $v_k(g)$ is assumed to be
equal to $\infty$ and $a_k(g)$ is not defined.
Let $\vv(g):= (v_1(g), \ldots, v_r(g))\in \Z^r_{\ge 0}$, $\aa(g):= (a_1(g), \ldots, a_r(g))\in (\C^*)^r$.

The functions (valuations) $v_k$ define a multi-index filtration on the ring
$\OO_{\C^n,0}$: for $\vv\in \Z^r$, the corresponding subspace $J(\vv)$ is defined
as
$\{g\in \OO_{\C^n,0} : \vv(g)\ge \vv\}$. In \cite{duke} there was computed the
(appropriately defined)
Poincar\'e series of the described multi-index filtration on the ring
$\OO_{\C^2,0}$ for a plane curve singularity $C$ (i.e., for $n=2$). It appeared to be equal to the Alexander polynomial of the
algebraic link $C\cap S^{3}_{\varepsilon}\subset S^3_{\varepsilon}$
corresponding to the curve $C$ ($S^3_{\varepsilon}$ is the sphere of radius
$\varepsilon$ centred at the origin of $\C^2$ with positive $\varepsilon$ small
enough).

Inspired by the notion of motivic integration (see, e.g., \cite{D-L},
\cite{loo})
there was defined the notion of the Euler
characteristic of (some) subsets of the ring  $\OO_{V,0}$  of functions on a
germ $(V,0)$ of an analytic variety or of its projectivization $\P\OO_{V,0}$ and
the corresponding notion of the integration with respect to the Euler
characteristic (see, e.g., \cite{stek2}). Here the Euler characteristic can be considered both as the usual one $\chi$ with values in $\Z$ and the generalized one $\chi_g$ with
values in the Grothendieck ring of complex algebraic varieties localized by the
class $\LL$ of the complex affine line: $K_0(\Var)_{(\LL)}$.

It was shown that the Poincar\'e series $P(\seq t1r)$ of a multi-index
filtration $\{J(\vv)\}$ on the ring $\OO_{V,0}$ (finitely determined: see the
definition below) is equal to the integral with respect to the Euler
characteristic over the projectivization $\P\OO_{V,0}$ of the ring of functions:
\begin{equation}\label{eq1}
P(\seq t1r) = \int_{\P\OO_{V,0}} \tt^{\,\vv(g)} d\chi
\end{equation}
($\tt=(t_1, \ldots, t_r)$, $\tt^{\,\vv} = t_1^{v_1}\cdot\ldots\cdot t_r^{v_r}$). Moreover, this representation
permitted to give a considerably shorter proof of the mentioned formula for the
Poincar\'e series of the filtration defined by orders $v_i(g)$ of a function on
the components of the plane curve singularity $C$ (see \cite{ijm}) and to compute
Poincar\'e series of some other multi-index filtrations:
\cite{div}, \cite{inv}, \dots

It is
natural to consider the integral like the one in the equation
(\ref{eq1}) with respect to the generalized Euler characteristic $\chi_g$.
This integral is a series in $\tt$ and $q = \LL^{-1}$. Generally speaking
it is a more fine invariant than the Poincar\'e series itself
(for non-plane curves). We give
two formulae  for it. The first one is in terms of the Hilbert
function $h(\vv)=\dim \OO_{V,0}/J(\vv)$. The second one, for the
described filtration defined by the orders of  functions on the components of
a plane curve singularity (and
also for the so called divisorial filtration for a modification
of $(\C^2,0)$ by a series of blowing-ups) is in terms of an
embedded resolution of the curve $C$ (or in terms of the
modification respectively).
We also give a formula for the generalized Euler characteristic of the extended semigroup
corresponding to the divisorial filtration.

\section{Multi-index filtrations and integrals with respect to
the Euler characteristic.}\label{sec1}

Let $(V,0)$ be a germ of an analytic variety. A one index
filtration
\begin{equation}\label{eq2}
\OO_{V,0} = J_0 \supset J_1 \supset \ldots \supset J_n\supset\ldots
\end{equation}
on the ring $\OO_{V,0}$ of germs of functions on $(V,0)$ by vector
subspaces can be defined by a function
$v:\OO_{V,0}\to \Z_{\ge 0}\cup \{\infty\}$
($v(g) = \sup \{i : g\in J_i\}$) with the properties
\begin{equation}
\label{eq3}
\begin{gathered}
v(\lambda\,g) = v(g) \quad \text{for } \lambda\in \C^*, \\
v(g_1+g_2) \ge \min \{v(g_1), v(g_2)\}.
\end{gathered}
\end{equation}
Sometimes the additional property
\begin{equation}
\label{eq4}
v(g_1g_2)=v(g_1)+v(g_2)
\end{equation}
is required (i.e., $v$ is a valuation). It guaranties that (\ref{eq2}) is a filtration by ideals and that $J_i\cdot J_j\subset J_{i+j}$.

The Poincar\'e series of the filtration (\ref{eq2}) is the
series
\begin{equation}
\label{eq5}
P(t) = \sum_{i = 0}^\infty \dim (J_i/J_{i+1})\cdot t^i
\end{equation}
(it makes sense if the dimensions of all the factors $J_i/J_{i+1}$ are
finite).

A multi-index filtration on the ring $\OO_{V,0}$ is defined by a
system $\{\seq v1r\}$ of functions
$\OO_{V,0}\to \Z_{\ge 0}\cup \{\infty\}$ satisfying the properties
(\ref{eq3}): the corresponding system of subspaces $J(\vv)$,
$\vv=(\seq v1r)\in \Z^r$, is defined by
$J(\vv) = \{g\in \OO_{V,0} : \vv(g)\ge \vv\}$. (Here
$\vv'\ge \vv''$ ($\vv', \vv''\in \Z^r$) if and only if
$v'_k\ge v''_k$ for all $k=1,\ldots, r$.)

\begin{remark}
It is sufficient to define the subspaces $J(\vv)$ only for
nonnegative $\vv$, i.e for $\vv\in\Z^r_{\ge 0}$, however, it will
be convenient if $J(\vv)$ is defined for all $\vv\in\Z^r$.
\end{remark}

We say that the filtration $\{J(\vv)\}$ is finitely determined if, for
any $\vv\in\Z^r$, there exists $N$ such that $J(\vv)\supset {\frak
m}^N$, where
${\frak m}$ is the maximal ideal in the ring $\OO_{V,0}$. If  the
multi-index filtration
$\{J(\vv)\}$ is finitely determined, all subspaces
$J(\vv)$ have finite codimension. In what follows all filtrations
are assumed to be finitely determined. Let $h(\vv) =\dim
\OO_{V,0}/J(\vv)$. We shall call $h(\vv)$ (as a function of $\vv$)
the Hilbert function of the filtration $\{J(\vv)\}$.

Before describing a multi-index version of the Poincar\'e
series it is convenient to define the notion of integration with respect
to the Euler characteristic over the projectivization
$\P\OO_{V,0}$ of the ring (space) $\OO_{V,0}$ and to interpret the
Poincar\'e series $P(t)$ as such an integral (similar to the
equation~(\ref{eq1})).

Let $K_0(\Var)$ be the Grothendieck ring of
quasi-projective varieties. It is generated by classes $[X]$ of
such varieties subject to the relations:\newline
1) if $X_1\cong X_2$, then $[X_1]=[X_2]$;\newline
2) if $Y$ is Zariski closed in $X$, then $[X]=[Y]+[X\setminus Y]$\newline
(the multiplication is defined by the Cartesian product).
Let $\LL$
be the class $[\A^1_{\C}]$ of the complex affine line.
The class $\LL$ is not equal to zero in the ring $K_0(\Var)$.
Moreover the natural ring homomorphism
$\Z[x]\to K_0(\Var)$ which sends $x$ to $\LL$ is an
inclusion.  Let
$K_0(\Var)_{(\LL)}$
be the localization of the
Grothendieck ring $K_0(\Var)$ by the class $\LL$.
The natural homomorphism
$\Z[x]_{(x)}\to K_0(\Var)_{(\LL)}$ is an inclusion as
well.

Let $J^p_{V,0}$, $p\ge 0$, be the space of $p$-jets of functions at the
origin in $(V,0)$, i.e. $J^p_{V,0}={\cal O}_{V,0}/{\frak
m}^{p+1}$, where $\frak m$ is the maximal
ideal in the ring ${\cal O}_{V,0}$.
It is a finite dimensional space. Let $d(p)=\dim J^p_{V,0}$. For
$(V,0)= (\C^n,0)$ one has $d(p) = \binom{n+p}{p}$.
For a complex vector space
$L$ (finite or infinite dimensional), let $\P
L=(L\setminus\{0\})/\C^*$ be its projectivization, let $\P^* L$ be
the disjoint union of $\P L$ with a point (in some sense $\P^*
L=L/\C^*$). One has natural maps $\pi_p: \P{\cal O}_{V,0} \to
\P^* J^p_{V,0}$ and $\pi_{p,q}: \P^* J^p_{V,0} \to \P^*
J^q_{V,0}$ for $p \ge q\ge 0$. Over $\P J^q_{V,0}
\subset \P^* J^q_{V,0}$ the map $\pi_{p,q}$ is a locally
trivial fibration, the fibre of which is a complex affine space
of dimension $d(p)-d(q)$.

\begin{definition}
A subset $X\subset \P{\cal O}_{V,0}$ is said to be cylindric if
$X=\pi_p^{-1}(Y)$ for a
constructible subset
$Y\subset \P J^p_{V,0} \subset \P^* J^p_{V,0}$.
\end{definition}

This definition means that the condition for a function
$g\in {\cal O}_{V,0}$ to belong to the set $X$ is a constructible
(and $\C^*$-invariant) condition on the $p$-jet of the
germ $g$.

\begin{definition}
For a cylindric subset $X\subset \P{\cal O}_{V,0}$
($X=\pi_k^{-1}(Y)$,
$Y\subset \P J^p_{V,0}$, $Y$ is constructible), its Euler
characteristic
$\chi(X)$
is defined as the Euler characteristic $\chi(Y)$ of the set $Y$.
The generalized Euler characteristic $\chi_g(X)$ of the cylindric
subset $X$ is the element $[Y]\cdot\LL^{-d(p)}$ of the  ring
$K_0(\Var)_{(\LL)}$.
\end{definition}

The generalized Euler characteristic $\chi_g(X)$ is well defined
since, if $X = \pi^{-1}_{q}(Y')$,
$Y'\subset \P J^{q}_{V,0}$, $p\ge q$, then $Y$ is  a
locally trivial fibration over $Y'$ with the fibre
$\C^{d(p)-d(q)}$ and therefore $[Y]=[Y']\cdot
\LL^{d(p)-d(q)}$.

\medskip
Let $\psi: \P{\cal O}_{V,0} \to A$ be a function with values in
an Abelian group $A$ with countably many values.

\begin{definition}
We say that the function $\psi$ is cylindric if, for each $a\ne 0$,
the set $\psi^{-1}(a)\subset \P{\cal O}_{V,0}$ is cylindric.
\end{definition}

\begin{definition}
The integral of a cylindric function $\psi$ over the space
$\P{\cal O}_{V,0}$ with respect to the Euler characteristic
(respectively with respect to the generalized Euler characteristic)
is
$$
\int_{\P{\cal O}_{V,0}}\psi d\chi_{\bullet} :=
\sum_{a\in A, a\ne 0} \chi_{\bullet}(\psi^{-1}(a))\cdot a,
$$
where $\chi_{\bullet}$ is the usual Euler characteristic $\chi$
or the generalized Euler characteristic $\chi_g$ respectively,
if this sum makes sense in
$A$
(respectively in
$A\otimes_{\Z}K_0(\Var)_{(\LL)}$). If the integral
exists
(makes sense) the function $\psi$ is said to be integrable.
\end{definition}

One can easily see that the Poincar\'e series (\ref{eq5}) of a
one index filtration is equal to the following integral with
respect to the Euler characteristic:
\begin{equation}\label{eq6}
P(t) = \int_{\P\OO_{V,0}} t^{\,v(g)}\, d\chi\; .
\end{equation}
Here $t^{\,v(g)}$ is considered as a function on $\OO_{V,0}$ (or
$\P\OO_{V,0}$) with values in the Abelian group $\Z[[t]]$ of
power series in $t$ with integral coefficients; $t^{\infty}$ is
supposed to be equal to $0$.
The above equality follows immediately from the fact that, for a
complex vector space $L$ of finite dimension, one has
$\,\dim L = \chi(\P L)$.

Let us use the following definition inspired by (\ref{eq6}).

\begin{definition}
The Poincar\'e series of a multi-index filtration $\{J(\vv)\}$,
$\vv\in\Z^r$, on the ring $\OO_{V,0}$ is the power series in $r$
variables $\seq t1r$ defined as the following integral with
respect to the Euler characteristic:
\begin{equation}\label{eq7}
P(\tt) = \int_{\P\OO_{V,0}} \tt^{\,\vv(g)}\, d \chi\; .
\end{equation}
\end{definition}
In \cite{duke} and \cite{div} it was essentially shown than the
Poincar\'e series $P(\tt)$ also has the following description
which in fact was its initial definition in \cite{CDK} (formally
speaking, in the cited papers this statement
is proved for a particular filtration, however, there is no
difference; see also Proposition~\ref{prop2} below).
Let $\1=(1, 1, \ldots, 1)\in\Z^r$.
Let
\begin{equation*}
L(\tt) = \sum_{\vv\in \Z^r} \dim(J(\vv)/J(\vv+\1))\cdot\tt^{\,\vv}\; .
\end{equation*}
The series $L(\tt)$ contains negative powers of variables $t_i$,
however, it does not contain monomials $\tt^{\,\vv}$ with purely
negative $\vv$ (i.e. all components of which are negative). It is
convenient to consider $L(\tt)$ as an element of the set
${\cal L}$ of formal Laurent series in the variables $\seq t1r$
with integer coefficients without purely negative exponents; i.e.
of expressions of the form
$\sum\limits_{\vv\in\Z^r\setminus \Z^r_{\le -1}}d(\vv)\cdot\tt^{\,\vv}$
with $d(\vv)\in\Z$.

\begin{proposition}\label{prop1}
\begin{equation}\label{eq8}
P(\tt) = \frac{L(\tt)\cdot\prod\limits_{k=1}^r(t_k-1)}{t_1\cdot\ldots\cdot t_r
-1}\,.
\end{equation}
\end{proposition}

In analogy to \cite{stek1}, \cite{div}, \cite{inv}, one can define a notion of
the extended semigroup of the multi-index filtration $\{J(\vv)\}$.
For $K\subset K_0=\{1, 2, \ldots, r\}$, let $\#K$ be the number of elements in
$K$, let $\1_K$ be the element of $\Z_{\ge0}^r$ whose $i$th component is equal
to $1$ or $0$ if $i\in K$ or $i\not\in K$ respectively (one has $\1=\1_{K_0}$).
For
$\vv\in\Z_{\ge0}^r$, let
$$
F_\vv:=(J(\vv)/J(\vv+\1))\setminus \bigcup_{k=1}^r
(J(\vv+\1_{\{k\}})/J(\vv+\1))\,.
$$
\begin{definition}
The extended semigroup $\widehat S$ is the union of the spaces $F_\vv$ for
$\vv\in\Z_{\ge0}^r$.
\end{definition}

The spaces $F_\vv$ are called fibres of the extended semigroup $\widehat S$.
The space $\widehat S$ is a graded space in the sense of \cite{duke} (i.e. its
components are numbered by elements of $\Z_{\ge0}^r$; in other words there is
a function $\vv$ on $\widehat S$ with values in $\Z_{\ge0}^r$ constant on
connected components $F_{\vv}$ of $\widehat{S}$: it has value
$\vv$ on $F_\vv$).
It is really a semigroup if the functions $v_i$  which define the filtration
are valuations, i.e. satisfy the property (\ref{eq4}). In this case the
semigroup operation in $\widehat S$ is defined by  multiplication of functions.
For the filtration defined by orders of functions on components $C_k$ of a curve $C= \bigcup\limits_{k=1}^r C_k\subset (\C^n, 0)$, the extended semigroup $\widehat S$ is the set $\{(\vv(g), \aa(g))\}\subset\Z^r_{\ge0}\times(\C^*)^r$ for all $g\in\O_{\C^n, 0}$ with $v_k(g)<\infty$ for $k=1,\ldots, r$.
The fibre $F_\vv$ is the complement to an arrangement of linear subspaces in a
linear space. Let $\P F_\vv=F_\vv/\C^*$ be its  projectivization. The union
$\P\widehat S$ of the spaces $\P F_\vv$ is (called) the projectivization of
the extended semigroup $\widehat S$. (The projectivization $\P\widehat S$ is a
semigroup itself if the functions $v_i$ are valuations.) From (\ref{eq7}) one can easily see that the Poincar\'e series $P(\tt)$ is equal to the integral of the function $\tt^{\,\vv}$ over the projectivization $\P\widehat S$ of the extended semigroup with respect to the Euler characteristic:
\begin{equation}\label{eq9}
P(\tt)=\int_{\P\widehat S} \ \tt^{\,\vv}\,d\chi.
\end{equation}

\section{Generalized Poincar\'e series of multi-index filtrations.}\label{sec2}

Using integration with respect to the generalized Euler characteristic one can
define the following two generalizations of the Poincar\'e series $P(t_1,\ldots,
t_r)$.

\begin{definition}
The {\it generalized Poincar\'e series} $P_g(\tt)$ of a filtration $J(\vv)$
defined by $\vv(g)=(v_1(g),\ldots,v_r(g))$ is the integral
$$
P_g(\tt)=\int_{\P\OO_{V,0}} \tt^{\,\vv(g)} d\chi_g\in
K_0(\Var)_{(\LL)}[[t_1,\ldots, t_r]]
$$
over the projectivization $\P\OO_{V,0}$ of the ring $\OO_{V,0}$
with respect to the {\it generalized} Euler characteristic $\chi_g$.
\end{definition}

The subset of the projectivization $\P\OO_{V,0}$ where
$\tt^{\,\vv(g)}$ is equal to $\tt^{\,\vv}$ (i.e., where $\vv(g)=\vv$)
is the projectivization of the complement
$J(\vv)\setminus\bigcup\limits_{k=1}^{r} J(\vv + \1_{\{k\}})$.
Because of that all the coefficients of the series
$P_g(\tt)$ are polynomials in $\LL^{-1}$. Therefore we shall write
$P_g(t_1,\ldots,t_r)$ as a series $P_g(t_1,\ldots,t_r,q)\in
\Z[[t_1,\ldots,t_r,q]]$ in $t_1$, \dots, $t_r$, and $q=\LL^{-1}$.
One has $P(t_1,\ldots,t_r)=P_g(t_1,\ldots,t_r,1)$.

To give ``a motivic version" of Proposition~\ref{prop1}, let us
define the corresponding version of the series $L(\tt)$:
$L_g(t_1,\ldots,t_r,q)\in{\cal L}[q]$. We have indicated that the
reason for Proposition~\ref{prop1} is the following formula for
the dimension. Let $A$ and $A'$ be subspaces of $\OO_{V,0}$ of
finite codimension. Then $\dim(A/A')=\chi(\P(A)\setminus\P(A'))$.
Substituting the usual Euler characteristic $\chi$ by the
generalized one $\chi_g$, one gets the following ``motivic version"
of the dimension: $``\dim_g"(A/A')=\chi_g(\P(A)\setminus\P(A'))$.
Let $\mbox{codim }A=a$, $\mbox{codim }A'=a'$ ($a' > a$). Then
$\chi_g(\P(A)\setminus\P(A'))=q^a+q^{a+1}+\ldots+q^{a'-1}=q^a\cdot\frac{1-q^{\,a'-a}}{1-q}$.
Therefore let $L_g(\tt,q)$ (``a motivic version" of $L(\tt)$) be
the series
$$
L_g(\tt,q):=\sum_{\vv\in\Z^r}
q^{\,h(\vv)}\cdot\frac{1-q^{\,h(\vv+\1)-h(\vv)}}{1-q}\cdot\tt^{\,\vv}\,.
$$

\begin{proposition}\label{prop2}
$$
P_g(\tt,q)=\frac{L_g(\tt,q)\cdot\prod\limits_{k=1}^r
(t_k-1)}{t_1\cdot\ldots\cdot t_r-1}\,.
$$
\end{proposition}

\begin{proof}
One has
\begin{eqnarray*}
P_g(\tt,q)&=&\int_{\P\OO_{V,0}}\tt^{\,\vv}\, d\chi_g
=\sum_{\vv\in\Z^r}\chi_g(\{g\in\P\OO_{V,0}:\vv(g)=\vv\})\cdot\tt^{\,\vv}=\\
&=&\sum_{\vv\in\Z^r}[\P J(\vv)\setminus\bigcup_{k=1}^r \P
J(\vv+\1_{\{k\}})]\cdot\tt^{\,\vv}=\\
&=&\sum_{{\vv}\in{\Z}^r}
\left(
\sum_{K\subset K_0}
(-1)^{{\#}K}
[\P J(\vv+\1_K)\setminus \P J(\vv+\1)]
\right)\cdot\tt^{\,\vv}\ .
\end{eqnarray*}

Therefore
\begin{multline*}
(t_1\cdot\ldots\cdot t_r-1)\cdot P_g(\tt, q)=\\
=\sum_{\vv\in\Z^r}
\left(
\sum_{K\subset K_0}(-1)^{\#K}
[\P J(\vv-\1+\1_K)\setminus \P J(\vv)]
\right)
\cdot\tt^{\,\vv}\ -\\
\qquad \qquad -\sum_{\vv\in\Z^r}
\left(
\sum_{K\subset K_0}(-1)^{\#K}
[\P J(\vv+\1_K)\setminus \P J(\vv+\1)]
\right)
\cdot\tt^{\,\vv}=\\
\ =\sum_{\vv\in\Z^r}
\sum_{K\subset K_0}(-1)^{\#K}[\P J(\vv-\1+\1_K)\setminus
\P J(\vv+\1_K)]\cdot\tt^{\,\vv}=\\
\qquad =\sum_{\vv\in\Z^r}\sum_{K\subset K_0}(-1)^{\#K}q^{\,h(\vv-\1+\1_K)}\cdot
\frac{1-q^{\,h(\vv+\1_K)-h(\vv-\1+\1_K)}}{1-q}\cdot\tt^{\,\vv}=\\
=\left(\sum_{\vv\in\Z^r}
q^{\,h(\vv)}\cdot\frac{1-q^{\,h(\vv+\1)-h(\vv)}}{1-q}\cdot\tt^{\,\vv}
\right
)\cdot\prod_{k=1}^r (t_k-1)=
L_g(\tt,q)\cdot\prod\limits_{k=1}^r(t_k-1)\ .
\end{multline*}
\end{proof}

Generally speaking, the generalized Poincar\'e series
$P_g(\tt,q)$ is a more fine invariant of a filtration than the
(usual) Poincar\'e series $P(\tt)$. Moreover this already holds
for filtrations defined by orders of functions on components of
curve singularities. Theorem~\ref{theo1} below shows that this
cannot happen for plane curve singularities: two plane curve
singularities with the same Poincar\'e series $P(\tt)$ are
topologically equivalent, have topologically equivalent
resolutions and therefore equal generalized Poincar\'e series
$P_g(\tt,q)$. However this can happen for non plane curve
singularities what can be seen from the following example taken from \cite{CDK}.

\begin{example}
Let $C= C_1\cup C_2\subset (\C^5,0)$ (respectively $C' =
C_1'\cup C'_2\subset (\C^6,0)$) be the germ of a curve whose
branches $C_1$ and $C_2$ (resp. $C'_1$ and $C'_2$) are defined by
the parametrizations ($t\in \C$, $u\in \C$):
\begin{eqnarray*}
C_1 &= \{(t^2,t^3,t^2,t^4,t^5)\}\,,\quad
C_2 &= \{(u^2,u^3, u^4, u^2, u^6)\}  \\
(\mbox{resp.}\;  C'_1 &=  \{(t^3,t^4,t^5, t^4,t^5, t^6)\}\,,\quad
C'_2 &= \{(u^3, u^4,  u^5, u^5, u^6, u^7)\} \; ).
\end{eqnarray*}
The local rings $\OO_{C,0}$ and $\OO_{C',0}$ of the curves $C$ and
$C'$, as subrings of
the normalizations $\C\{t\}\times \C\{u\}$, are:
\begin{eqnarray*}
\OO_{C,0} & = & \C\{\,(t^2,u^2), (t^3, u^3), (t^2,u^4), (t^4, u^2), (t^5,
u^6)\,\}, \\
\OO_{C',0} & = & \C\{\,(t^3,u^3), (t^4,u^4), (t^5, u^5), (t^4,u^5),
(t^5, u^6), (t^6, u^7)\,\}\; .
\end{eqnarray*}
The corresponding semigroups of values are, respectively,
\begin{eqnarray*}
S &= &\{(0,0), (3,3), (2,k), (\ell,2), (r,s) : k\ge 2, \ell\ge 3,
r, s \ge 4\}\,, \\
S\,' &= &\{(0,0), (3,3), (r,s) : r, s \ge 4\} \; .
\end{eqnarray*}
An easy computation shows that in both cases the Poincar\'e
series is the polynomial $P(\tt) = 1 + t_1^{3}t_2^3$. However
one has $h(3,3)=1$ for the curve $C'$ and $h(3,3)=3$ for $C$
(note that a basis for this case consists of the classes of
$1$, $(t^2,u^4)$ and $(t^4,u^2)$ with the linear
dependence
$(t^2,u^4)+(t^4,u^2)-(t^2,u^2)=(t^4,u^4)$ in $J(3,3)$ ).
\end{example}

\begin{remark}
The example shows that, generally speaking, the Poincar\'e series $P(t_1,\ldots, t_r)$ of the multi-index filtration $\{J(\vv)\}$ defined by a system of functions $\{v_k\}$, $k\in K_0=\{1,2,\ldots, r\}$, does not determine the Hilbert function $h(\vv)$. However, the set of the Poincar\'e series of the filtrations defined by all subsystems of $\{v_k\}$ does determine the Hilbert function. For $K\subset K_0$, let $P_K(\tt)$ be the Poincar\'e series corresponding to the system $\{v_k:\, k\in K\}$ (it really depends on variables $t_k$ only with $k\in K$; $P_{K_0}(\tt)=P(\tt)$). Let $H(\tt)=\sum\limits_{\vv\in\Z_{\ge 0}^r}\cdot h(\vv)\tt^\vv\,$. The same arguments as in \cite{div} or \cite{inv} show that
$$
H(\tt)=\frac{\tt}{1-\tt^{\1}}\cdot
\frac
{\sum\limits_{K\subset K_0}(-1)^{r-\# K}P_K(\tt)\cdot(t_1\cdot\ldots\cdot t_r-1)
\mbox{
\raisebox{-0.5ex}{$\vert$}}{}_{\{t_k=1 \mbox{ for } k\not\in K\}}
}{\prod_{k=1}^r(t_k-1)}
$$
(in \cite{div} there is a misprint: the factor $\tt$ is forgotten; it is corrected in \cite{inv}).
The fact that, in general, the Poincar\'e series $P(\tt)$ does not determine the Hilbert function $h(\tt)$ implies that it does not determine the Poincar\'e series $P_K(\tt)$ for $K\subset K_0$. This can be explained by the fact that both $P(\tt)$ and $P_K(\tt)$ can be computed as integrals of the type (\ref{eq7}). However, those $g\in\O_{V,0}$, for which $v_k(g)<\infty$ for all $k \in K$ but $v_k(g)=\infty$ for some $k\not\in K$, participate in the integral for $P_K(\tt)$, but do not participate in the integral for $P(\tt)$: in this case $\tt^{\,\vv(g)}=0$. This gives the following fact. Suppose that the set of functions $v_k$ which defines the filtration is such that if $v_k(g)=\infty$ for certain $k\in K_0$ then $v_k(g)=\infty$ for all $k\in K_0$. Then $P_K(\tt)=P(\tt)\mbox{
\raisebox{-0.5ex}{$\vert$}}{}_{\{t_k=1 \mbox{ for } k\not\in K\}}$.
This holds, in particular, for the divisorial valuations discussed below.
\end{remark}

One can also consider a generalization of the Poincar\'e series
$P(t_1,\ldots,t_r)$ which emerges from the equation~(\ref{eq9}).

\begin{definition}
The {\it generalized semigroup Poincar\'e series} ${\widehat P}_g(\tt)$ of a
filtration $\{J(\vv)\}$ defined by $\vv(g)=(v_1(g),\ldots,v_r(g))$ is the
integral
$$
{\widehat P}_g(\tt)=\int_{\P\widehat S} \tt^{\,\vv(g)}\,d\chi_g\in
K_0(\Var)[[t_1,\ldots, t_r]]
$$
over the projectivization of the extended semigroup $\widehat S$
with respect to the {\it generalized} Euler characteristic $\chi_g$.
\end{definition}

All connected components of $\P\widehat S$ (i.e. projectivizations $\P F_\vv$ of the fibres $F_\vv$\,) are complements to arrangements of projective subspaces in {\it finite dimensional} (!) projective spaces. Because of that all the coefficients of the series
${\widehat P}_g(\tt)$ are polynomials in $\LL$. Therefore we shall write
${\widehat P}_g(t_1,\ldots,t_r)$ as a series ${\widehat P}_g(t_1,\ldots,t_r,\LL)\in
\Z[[t_1,\ldots,t_r,\LL]]$ in $t_1$, \dots, $t_r$, and $\LL$.
One has $P(t_1,\ldots,t_r)={\widehat P}_g(t_1,\ldots,t_r,1)$.

Repeating the arguments of Proposition~\ref{prop2}, one gets the following equation:
$$
\widehat P_g(\tt)=\frac{\left(\sum\limits_{\vv\in\Z^r}[\P(J(\vv)/J(\vv+1))]\cdot\tt^{\,\vv}\right)\cdot
\prod\limits_{k=1}^r(t_k-1)}{t_1\cdot\ldots\cdot t_r -1}\ .
$$

\section{The Poincar\'e series $P_g(\tt,q)$ in terms of an embedded
resolution.}\label{sec3}

Let $C\subset(\C^2,0)$ be a (reduced) plane curve singularity,
$C=\bigcup\limits_{k=1}^r C_k$,
$C_k=\{f_k=0\}$, $f_k\in\OO_{\C^2,0}$, and
let $\pi:(\X,\D)\to(\C^2,0)$ be an embedded resolution of it. Let
the exceptional divisor $\D$ of the resolution $\pi$ be the union
of irreducible components $E_i$, $i=1,\ldots,s$, each of them is
isomorphic to the projective line $\C\P^1$. Let
$\stackrel{\circ}{E_i}$ be the ``nonsingular part" of the component
$E_i$, i.e., $E_i$ minus intersection points with all other
components of the total transform $(f\circ\pi)^{-1}(0)$ of the
curve $C$. We shall denote by $\stackrel{\bullet}{E_i}$ the
``nonsingular part" of the component $E_i$ in the space $\X$
of the resolution, i.e., $E_i$ minus intersection points with all
other components $E_{j}$ of the exceptional divisor $\D$. For
$i=1,\ldots,s$ and $g\in\OO_{\C^2,0}\setminus\{0\}$, let $w_i(g)$
be the multiplicity of the lifting $\widetilde g=g\circ\pi$ of the
function $g$ to the space $\X$ of the resolution along the
component $E_i$ of the exceptional divisor $\D$.
The valuations $w_i$ define a multi-index filtration (called divisorial)
$\{J^D(\ww)\}$ on the ring
$\OO_{\C^2,0}$: $J^D(\ww) =\{g :\, \ww(g)\ge \ww\}$.
For
$k=1,\ldots,r$, let $i(k)$ be the number of the component $E_i$
($i=i(k)$) of the exceptional divisor $\D$ which intersects the
strict transform $\widetilde C_k$ of the component $C_k$ of the
curve $C$. Let $A=(E_i\circ E_j)$ be the intersection matrix of
the components $E_i$. (For each $i$ the self-intersection number
$(E_i\circ E_i)$ is negative, for $i\ne j$ the intersection number
$(E_i\circ E_j)$ is either $0$ or $1$; $\det A=\pm 1$.) Let
$M=-A^{-1}$. The entries $m_{ij}$ of the matrix $M$ are positive
and have the following meaning. Let $\widetilde L_i$ be a germ of
a smooth curve on $\X$ transversal to the component $E_i$ of the
exceptional divisor $\D$ at a smooth point (i.e., at a point of
$\stackrel{\bullet}{E_j}$) and let the projection
$L_i=\pi(\widetilde L_i)\subset(\C^2,0)$ of the curve $\widetilde
L_i$ be given by an equation $g_i=0$, $g_i\in\OO_{\C^2,0}$. Then
$m_{ij}=w_j(g_i)=w_i(g_j)=L_i\circ L_j$.

Let $I_0=\{(i,j): i<j,\, E_i\cap E_j=pt\}$, $K_0=\{1, 2,\ldots,r\}$.
For $\sigma\in I_0$, $\sigma=(i,j)$, let $P_\sigma=E_i\cap E_j$;
for $k\in K_0$, let $P_k=E_{i(k)}\cap\widetilde C_k$. For $\sigma=(i,j)\in
I_0$,
let $i(\sigma):=i$, $j(\sigma):=j$. For $I\subset I_0$, $K\subset K_0$, let
\begin{eqnarray*}
{\cal N}_{I, K}&:=&\{\n=(n_i, n'_\sigma, n''_\sigma, {\widetilde n}'_k,
{\widetilde n}''_k): n_i\ge 0, i=1,\ldots, s;\\&&
n'_\sigma>0, n''_\sigma>0, \sigma\in I;\
{\widetilde n}'_k>0, {\widetilde n}''_k>0, k\in K\}.
\end{eqnarray*}
For $\n\in{\cal N}_{I, K}$, $i=1, \ldots, s$, let
$$
{\widehat n}_i := n_i + \sum_{\sigma\in I:\,i(\sigma)=i}n'_\sigma +
\sum_{\sigma\in I:\, j(\sigma)=i}n''_\sigma +
\sum_{k\in K:\, i(k)=i}{\widetilde n}'_k.
$$
Let
\begin{equation}
\label{eq}
F(\n) := \frac{1}{2}\left(
\sum_{i,j=1}^s m_{ij}{\widehat n}_i {\widehat n}_j + \sum_{i=1}^s {\widehat n}_i\cdot
\left(
\sum_{j=1}^s m_{ij}\chi(\stackrel{\bullet}{E_j}) + 1
\right)
\right)
+\sum_{k\in K}{\widetilde n}''_k\ ,
\end{equation}
$\ww(\n):= \sum\limits_{i=1}^s {\widehat n_i}\,\mm_i$, where $\mm_i=(m_{i1}, \ldots,
m_{is}) \in\Z_{\ge0}^s$, $v_k(\n):= w_{i(k)}(\n)+{\widetilde n}''_k$.

\begin{theorem}\label{theo1}
\begin{eqnarray*}
P_g(t_1,\ldots,t_r,q)&=&\sum_{I\subset I_0,{\ }
K \subset K_0
}
\quad
\sum_{\n\in{\cal N}_{I, K}} q^{F(\n)-\sum\limits_{i=1}^s
n_i-\#I-\#K}
\cdot(1-q)^{\#I+\#K}\times
\\
&\times&
\prod_{i=1}^s\left(\sum\limits_{j=0}^{\min\{n_i,
1-\chi(\stackrel{\circ}{E_i})\}}
(-1)^j\binom{1-\chi(\stackrel{\circ}{E_i})}{j}q^j \right)\cdot\tt^{\,\vv(\n)}.
\end{eqnarray*}
\end{theorem}

\begin{remark}
One can see that the generalized Poincar\'e series $P_g(\tt, q)$ of the discussed filtration represents a rational function in the variables $t_1$, \dots, $t_r$, and $q$ (in fact with the denominator $\prod\limits_{k=1}^r (1-qt_k)$). This is not immediately clear from the equation of Theorem~\ref{theo1}.
\end{remark}

\begin{proof}
In a neighbourhood of the point $P_\sigma$, $\sigma\in I_0$, (respectively of
the point $P_k$, $k\in K_0$) choose local coordinates $x_\sigma$, $y_\sigma$
(respectively $x_k$, $y_k$) so that the components of the total transform
$\pi^{-1}(C)$ of the curve $C$ are the coordinate lines:
$E_{i(\sigma)}=\{y_\sigma=0\}$, $E_{j(\sigma)}=\{x_\sigma=0\}$
(respectively $E_{i(k)}=\{y_k=0\}$, ${\widetilde C}_k=\{x_k=0\}$).

For a space $X$, let $S^n X=X^n/S_n$ be the $n$th symmetric power of the space $X$.
Let
$$
Y:=\bigcup_{I\subset I_0,{\ }K\subset K_0}
\quad
\bigcup_{\n\in{\cal N}_{I, K}}
\left(
\prod_{i=1}^s S^{n_i} \stackrel{\circ}{E_i}\times \prod_{\sigma\in I}\C^*_\sigma
\times \prod_{k\in K}\C^*_k
\right),
$$
where
$\C^*_\sigma$ and $\C^*_k$ are copies of the punctured line $\C^*=\C\setminus\{0\}$ numbered by $\sigma\in I_0$ and $k\in K_0$
respectively. Let $Y_\n=\prod_{i=1}^s S^{n_i} \stackrel{\circ}{E_i}\times
\prod_{\sigma\in I}\C^*_\sigma \times
\prod_{k\in K}\C^*_k$ be the corresponding connected component of the space
$Y$. The space $Y$ can be considered as a semigroup with respect to the
following semigroup operation. Let $y_1$ and $y_2$ be two points of $Y$ which
belong to the components $Y_{{}^1{\n}}$ and $Y_{{}^2{\n}}$ respectively, where
${}^1{\n}=({}^1 n_i, {}^1n'_\sigma, {}^1 n''_\sigma, {}^1{\widetilde n}'_k,
{}^1{\widetilde n}''_k)\in {\cal N}_{I^1, K^1}$,
${}^2 \n=({}^2 n_i, {}^2 n'_\sigma, {}^2 n''_\sigma, {}^2{\widetilde n}'_k,
{}^2 {\widetilde n}''_k)\in {\cal N}_{I^2, K^2}$.
Then the product $y_1y_2$ belongs to the component $Y_{{\n}}$ with $\n=(n_i,
n'_\sigma, n''_\sigma, {\widetilde n}'_k,{\widetilde n}''_k)\in{\cal N}_{I,
K}$ such that $I=I^1\cup I^2$, $K=K^1\cup K^2$, each of $n_i$, $n'_\sigma$,
$n''_\sigma$, ${\widetilde n}'_k$, and ${\widetilde n}''_k$ is equal to the
sum of the corresponding components in ${}^1 \n$ and ${}^2 \n$ where, if a
component is absent (say, if $\sigma\in I$, but $\sigma\not\in I^1$), it is
assumed to be equal to zero. Moreover, in the factor $S^{n_i}
\stackrel{\circ}{E_i}$ of the component $Y_\n$ the product $y_1y_2$ is
represented by the union of the corresponding tuples of points of
$\stackrel{\circ}{E_i}$ in $y_1$ and in $y_2$ (if a component is absent, the
corresponding tuple is assumed to be empty); in the factor $\C^*_\sigma$ or
$\C^*_k$ the product $y_1y_2$ is represented by the product of the
corresponding components (nonzero numbers) in $y_1$ and in $y_2$ (if a
component is absent, it is assumed to be equal to $1$).

Let $\OO^*_{\C^2,0}=\{g\in\OO_{\C^2,0}: v_i(g)<\infty\mbox{ for }i=1,\ldots,r\}$
and let $\Init$ be the map $\P\OO^*_{\C^2,0}\to Y$ from the projectivization of
$\OO^*_{\C^2,0}$ defined in the following way. For $g\in \OO^*_{\C^2,0}$, let
$\Gamma_g\subset\X$ be the strict transform of the curve
$\{g=0\}\subset(\C^2,0)$. Let $I(g):=\{\sigma\in I_0: P_\sigma\in\Gamma_g\}$,
$K(g):=\{k\in K_0: P_k\in\Gamma_g\}$. Let $n_i=n_i(g)$, $i=1,\ldots,s$,
be the number of intersection points of the curve $\Gamma_g$ with the smooth
part
$\stackrel{\circ}{E_i}$ of the component $E_i$ counted with multiplicities.
For $\sigma\in I(g)$, let $n'_\sigma=n'_\sigma(g):=(\Gamma_g\circ
E_{i(\sigma)})_{P_\sigma}$,
$n''_\sigma=n''_\sigma(g)=(\Gamma_g\circ E_{j(\sigma)})_{P_\sigma}$;
for $k\in K(g)$, let ${\widetilde n}'_k={\widetilde n}'_k(g):=
(\Gamma_g\circ E_{i(k)})_{P_k}$, ${\widetilde n}''_k={\widetilde n}''_k(g):=
(\Gamma_g\circ {\widetilde C}_k)_{P_k}$.
For $\sigma\in I(g)$ (respectively for $k\in K(g)$),
in a neighbourhood of the point $P_\sigma$ (respectively of the point $P_k$),
the germ of the curve $\Gamma_g$ is given by an equation
$\varphi_\sigma(x_\sigma,y_\sigma)=0$, where
$\varphi_\sigma(x_\sigma,y_\sigma)=
a_\sigma(g)x_\sigma^{n'_\sigma}+b_\sigma(g)y_\sigma^{n''_\sigma}+\ \mbox{terms
with monomials } x_\sigma^\alpha y_\sigma^\beta \mbox{ with either } \alpha>0,
\beta>0,\mbox{ or } \alpha>n'_\sigma,\mbox { or }\beta>n''_\sigma$,
$a_\sigma(g)\ne 0$, $b_\sigma(g)\ne 0$
(respectively by an equation $\varphi_k(x_k,y_k)=0$, where
$\varphi_k(x_k,y_k)=
a_k(g)x_k^{{\widetilde n}'_k}+
b_k(g)y_k^{{\widetilde n}''_k}+\ldots$\,, $a_k(g)\ne 0$, $b_k(g)\ne 0$). Now
$\Init(g)$ is an element of the summand $Y_\n$ corresponding to $I(g)$, $K(g)$
and $\n=(n_i(g), n'_\sigma(g), n''_\sigma(g), {\widetilde n}'_k(g),
{\widetilde n}''_k(g))$ the components (factors) of which are the following
ones. In $S^{n_i} \stackrel{\circ}{E_i}$, $i=1,\ldots, s$, it is represented
by the set of the intersection points of the curve $\Gamma_g$ and
$\stackrel{\circ}{E_i}$ counted with multiplicities; in $\C^*_\sigma$,
$\sigma\in I(g)$, (respectively in $\C^*_k$, $k\in K(g)$)
it is represented by the element $a_\sigma(g):b_\sigma(g)$ (respectively by
$a_k(g):b_k(g)$). One can see that the map $\Init$ is a semigroup homomorphism
$\P\OO^*_{\C^2,0}\to Y$, where the semigroup structure on $\P\OO^*_{\C^2,0}$ is
defined by multiplication of functions.

\begin{lemma}\label{lemma1}
For any point $y\in Y$, the preimage $\Init^{-1}(y)$ is an affine space in the projectivization $\P\O_{\C^2,0}$ of a finite codimension. Moreover, over each connected component $Y_\n$ of the space $Y$ the map $Int$ is a locally trivial fibration.
\end{lemma}

This statement follows from the next one.

\begin{lemma}\label{lemma2}
For $g$ and $g'$ from $\OO^*_{\C^2,0}$, $\Init(g)=\Init(g')$ if and only if
$g'=\alpha g+h$, where $\alpha\in\C^*$, $w_i(h)>w_i(g)(=w_i(g'))$ for all
$i=1,\ldots,s$,
$v_k(h)>v_k(g)(=v_k(g'))$ for all $k=1,\ldots,r$.
\end{lemma}

\begin{proof}
Let $\widetilde g=g\circ\pi$ and ${\widetilde g}{\,}'=g'\circ\pi$ be the liftings of
the functions $g$ and $g'$ to the space $\X$ of the resolution.  Their ratio
$\Psi=\widetilde g{\,}'/\widetilde g$ is a meromorphic function on $\X$. If
$\Init(g)=\Init(g')$, then zeroes and poles of the function $\Psi$ on each
component of the exceptional divisor $\D$ and of the strict transform
$\widetilde C$ of the curve $C$ cancel each other and therefore $\Psi$ is a
regular function on each of them (in particular, it is
constant on each component $E_i$). The conditions on the coefficients in the
power series
decompositions of the functions $\widetilde{g}$ and ${\widetilde{g}}{\,}'$ at
the intersection
points of the total transform $\pi^{-1}(C)$ of the curve $C$ imply that at all
these points the values of the function $\Psi$ on the both components of the
total transform coincide. Therefore $\Psi$ is a regular function on
$\pi^{-1}(C)$ equal to a constant (say, $\alpha$) on the exceptional divisor
$\D$. Therefore $w_i(g'-\alpha g)>w_i(g)\,(\,=w_i(g')\,)$,
$i=1,\ldots,s$, $v_k(g'-\alpha g)>v_k(g)\,(\,=v_k(g')\,)$, $k=1,\ldots,r$.

In the other direction, if $g'=\alpha g+h$ with $\alpha\in\C^*$,
$w_i(g'-\alpha g)>w_i(g)$, $i=1,\ldots,s$, $v_k(g'-\alpha g)>v_k(g)$,
$k=1,\ldots,r$, then $\Psi$ is a regular function on the total transform
$\pi^{-1}(C)$ and $\Psi_{\vert\D}\equiv\alpha$.
Therefore zeroes and poles of the function $\Psi$ cancel each other on each
component of the exceptional divisor $\D$ and therefore the intersection
points of the strict transform of the curves $\{g=0\}$ and $\{g'=0\}$ with
each component $E_i$ of the exceptional divisor $\D$ coincide (with their
multiplicities). The intersection number of the strict transform of the curve
$\{g=0\}$ with the component $\widetilde C_k$ of the strict transform of the
curve $C$ is equal to $v_k(g)-w_{i(k)}(g)$ and therefore coincides with that
of $\{g'=0\}$.
The fact that the ratios of the coefficients in local equations of the strict
transforms of the curves $\{g=0\}$ and $\{g'=0\}$ at the intersection points
of the components of $\pi^{-1}(C)$ are equal is obvious.
\end{proof}

The number $w_i(\n)$ defined before Theorem \ref{theo1} is just the multiplicity along
the component $E_i$ of the exceptional divisor ${\cal D}$ of the lifting of a
function $g$ such that $\Init(g)\in Y_{\n}$. The number $v_k(\n)$ is the order
of such a function on the component $C_k$ of the curve.

Lemma \ref{lemma2} says that the preimage $\Init^{-1}(y)$ of a point
$y\in Y_{\n}$ is an affine space whose codimension $F(\n)$ in $\P\OO_{\C^2,0}$ is one less
than the codimension of the ideal $I_{\n}$ in $\OO_{\C^2,0}$ which
consists of functions $h$ with $w_i(h)>w_i(\n)$, $i=1,\ldots, s$ and
$v_k(h) >v_k(\n)$, $k=1,\ldots,r$. Applying the Fubini formula to the map $\Init$
one gets
$$
P_g(\tt,q) = \int_{Y} q^{\,F(\n)} \tt^{\,\vv(\n)} d\chi_g =
\sum_{I\subset I_0,{\ }K\subset K_0}\quad
\sum_{\n\in{\cal N}_{I, K}} q^{F(\n)} \; [Y_{\n}] \; \tt^{\vv(\n)}\; .
$$
Here $[Y_{\n}]=[\C^*]^{\# I +\# K}\cdot \prod_{i=1}^s
[S^{n_i}\stackrel{\circ}{E_i}]$,
$[\C^*] = q^{-1}-1 = q^{-1}(1-q)$. The class
$[S^{n_i}\stackrel{\circ}{E_i}]$ as a polynomial in $\LL=q^{-1}$ can be
obtained from the formula
\begin{eqnarray*}
\sum_{n=0}^{\infty} [S^n\stackrel{\circ}{E_i}]\, t^n &=& (1-t)^{-[\LL -
(1-\chi(\stackrel{\circ}{E_i}))]} \\
&=& (1-t)^{-\LL}\cdot (1-t)^{1-\chi(\stackrel{\circ}{E_i})} =
\sum_{n=0}^{\infty} \LL^{n} t^n \cdot
(1-t)^{1-\chi(\stackrel{\circ}{E_i})}
\end{eqnarray*}
(see, e.g., \cite{glm}). Therefore
\begin{eqnarray*}
[S^{n}\stackrel{\circ}{E_i}] &=&
\sum_{j=0}^{\min(n,1-\chi(\stackrel{\circ}{E_i}))}
(-1)^j \LL^{n-j} \binom{1-\chi(\stackrel{\circ}{E_i})}{j} \\
\ &=&
q^{-n}\sum_{j=0}^{\min(n,1-\chi(\stackrel{\circ}{E_i}))}
(-1)^j  \binom{1-\chi(\stackrel{\circ}{E_i})}{j} q^j\; .
\end{eqnarray*}

Now to prove Theorem \ref{theo1} it remains to show that the codimension
$F(\n)$ of the ideals $I_{\n}$ is given by the formula (\ref{eq}). Making
additional blowing-ups at intersection points of the strict transforms
$\widetilde{C_k}$ of the curve $C$ with the exceptional divisor ${\cal D}$ one
reduces this problem to the case when $K=\emptyset$. In this case the
codimension we are interested in is equal to $h^D(\ww(\n)+\1)-1$ where
$h^D(\ww) = \dim \OO_{\C^2,0}/J^D(\ww)$ is the Hilbert function of the
divisorial filtration $\{J^D(\ww)\}$.
The codimension $h^D(\ww(\n))$ can be computed by the formula
$h^D(\ww(\n)) = -\frac 12 (\DD \circ \DD + \DD \circ \K)$ (it can be also obtained from
the Hoskin--Deligne formula: see, e.g, \cite{casas}), where
$\DD = \sum_{i=1}^s w_i(\n)\, E_i= - \sum_{i=1}^s \widehat{n}_i\, E^{*}_i$,
$\K = - \sum_{i=1}^s (2+ E_i\circ E_i) E^{*}_i$ is the canonical divisor of $\X$
and $\seq{E^*}1s$ is the basis of divisors on $\X$ dual (with respect to the
intersection form) to $\seq E1s$.
Finally the formula for $h^D(\ww(\n)+\1)$ follows from the fact that an open
part of $\P J^D(\ww(\n))$ is fibred over the $\sum_{i=1}^s \widehat{n}_i$
dimensional space $\prod_{i=1}^s S^{\widehat{n}_i}\stackrel{\bullet}{E_i}$ with
the fibre $J^D(\ww(\n)+\1)$ (Lemma~\ref{lemma2} applied to a function $g$ with
$I(g)=\emptyset$).

\begin{remark}
If $\ww(\n)+\1$ itself belongs to the semigroup of values of the set of divisorial valuations $\{\ww_k\}$ (i.e., if $\ww(\n)+\1=\ww(\n'$), the last argument can be substituted by the direct computation of $h^D(\ww(\n)+\1)$ with the use of the formula indicated above.
\end{remark}
\end{proof}

\section{The Poincar\'e series $P_g(\tt,q)$ for divisorial
valuations and the Euler characteristic of its extended semigroup.}
\label{sec4}

Let $\pi:(\X,\D)\to(\C^2,0)$ be a modification of the complex plane
by a sequence of blowing-ups at preimages of the origin in $\C^2$.
Let the exceptional divisor $\D$ of the modification $\pi$ be the union
of irreducible components $E_i$, $i=1,\ldots,s$, each of them is
isomorphic to the projective line $\C\P^1$. As above, let
$\stackrel{\bullet}{E_i}$ be the ``nonsingular part" of the component
$E_i$, i.e., $E_i$ minus intersection points with all
other components $E_{j}$ of the exceptional divisor ${\D}$. For
$i=1,\ldots,s$ and $g\in\OO_{\C^2,0}\setminus\{0\}$,
let $w_i(g)$
be the multiplicity of the lifting ${\widetilde g}=g\circ\pi$ of the
function $g$ to the space $\X$ of the resolution along the
component $E_i$ of the exceptional divisor $\D$. Let $A=(E_i\circ E_j)$
be the intersection matrix of the components $E_i$ and let
$M= (m_{ij})= -A^{-1}$. Let $\{J^D(\ww)\}$ be the divisorial filtration defined
by the valuations $\ww_i(g)$, $i=1, \ldots, s$: $J^D(\ww):=\{g\in\OO_{\C^2,0}:
\ww(g)\ge\ww\}$.

Let $I_0=\{(i,j):\, i<j,\, E_i\cap E_j=pt\}$. For $\sigma\in I_0$, $\sigma=(i,j)$,
let $P_\sigma$, $i(\sigma)$, and $j(\sigma)$ be defined as in section~\ref{sec3}.
For $I\subset I_0$, let
\begin{eqnarray*}
{\cal N}^D_{I}:=\{\n=(n_i,\, n'_\sigma,\, n''_\sigma :\, n_i\ge 0,\, i=1,\ldots,
s;\ n'_\sigma>0,\, n''_\sigma>0,\, \sigma\in I\}.
\end{eqnarray*}
For $\n\in{\cal N}^D_{I}$, $i=1, \ldots, s$, let
$$
{\widehat n}_i := n_i + \sum_{\sigma\in I:\, i(\sigma)=i}n'_\sigma +
\sum_{\sigma\in I:\, j(\sigma)=i}n''_\sigma.
$$
Let
$$
F^D(\n) := \frac{1}{2}\left(
\sum_{i,j=1}^s m_{ij}{\widehat n}_i {\widehat n}_j + \sum_{i=1}^s {\widehat n}_i\cdot
\left(
\sum_{j=1}^s m_{ij}\chi(\stackrel{\bullet}{E_j}) + 1
\right)
\right),
$$
$\ww(n):= \sum\limits_{i=1}^s {\widehat n}_i\, \mm_{\,i}$.

\begin{theorem}\label{theo2}
The generalized Poincar\'e series $P^D_g(t_1,\ldots,t_r,q)$ of the divisorial
filtration
$\{J^D(\ww)\}$ is equal to
\begin{eqnarray*}
P^D_g(t_1,\ldots,t_r,q)
&=&\sum_{I\subset I_0}
\quad
\sum_{\n\in{\cal N}^D_{I}} q^{F^D(\n)-\sum\limits_{i=1}^s n_i-\#I}
\cdot(1-q)^{\#I}\times
\\
&\times&\prod_{i=1}^s \left(\sum\limits_{j=0}^{\min\{n_i, 1-\chi(\stackrel{\bullet}{E_i})\}} (-1)^j\binom{1-\chi(\stackrel{\bullet}{E_i})}{j}q^j \right)
\cdot\tt^{\,\ww(\n)}\ .
\end{eqnarray*}
\end{theorem}

\begin{proof}
The arguments essentially repeat those of Theorem~\ref{theo1} with the use of
the space
$$
Y^D:=\bigcup_{I\subset I_0}
\quad
\bigcup_{\n\in{\cal N}^D_{I}}
\left(
\prod_{i=1}^s S^{n_i} \stackrel{\bullet}{E_i}\times \prod_{\sigma\in I}\C^*_\sigma
\right).
$$
\end{proof}

In \cite{div} it was shown that
\begin{equation}\label{eq10}
P^D(t_1,\ldots, t_s)=\chi(\P \widehat
S_D)=\prod_{i=1}^s\left(1-
\tt^{\,\mm_{\,i}}\right)^{-\chi(\stackrel{\bullet}{E_i})}
\end{equation}
(it follows easily from the formula of Theorem~\ref{theo2}). Since
$\chi(\stackrel{\bullet}{E_i})=2-\#\{\sigma\in I_0: i(\sigma)=i \mbox{ or }
j(\sigma)=i\}$, the equation~\ref{eq10} can be rewritten as
\begin{equation}\label{eq11}
\chi(\P \widehat S_D)=\frac{\prod_{\sigma\in
I_0}(1-\tt^{\,\mm_{\,i(\sigma)}})(1-\tt^{\,\mm_{\,j(\sigma)}})}
{\prod_{i=1}^s (1-\tt^{\,\mm_{\,i}})^2}\ .
\end{equation}
Let us write a version of the equation~(\ref{eq11}) for the generalized Euler characteristic.

Using the construction described in \cite{duke} one can define a map (in fact
a semigroup homomorphism) $\Pi:Y^D\to\P {\widehat S}_D$ in the following way. For
$y\in Y$, let
$\Gamma$ be an arbitrary (germ of a) curve on the space $(\X,\D)$ of the
resolution which represents the point $y$: see the description of the point of
the space $Y$ corresponding to the curve $\Gamma_g$ in the proof of
Theorem~\ref{theo1}. Let the projection $\pi(\Gamma)$ of the curve $\Gamma$ be
given by an equation $\{g=0\}$. If there is another curve $\Gamma'$ like that
and $\pi(\Gamma')=\{g'=0\}$, then the meromorphic function
$\Psi=g'\circ\pi / g\circ\pi$ on the space $\X$ of the resolution is
constant on the exceptional divisor $\D$. This implies that the points of the
projectivization $\P {\widehat S}_D$ of the extended semigroup corresponding to
the functions $g$ and $g'$ coincide. By definition this point is the image
$\Pi(y)$ of the point $y$.

\begin{proposition} \label{prop5}
$\Pi$ is a semigroup isomorphism.
\end{proposition}

\begin{proof}
This is a direct consequence of Lemma~\ref{lemma2}.
\end{proof}

\begin{theorem} \label{theo3}
$$
{\widehat P}^D_g(\tt,\LL)=\chi_g(\P \widehat S_D) = \int_{\P \widehat S_D}
\tt^{\,\ww}\ d \chi_{g} =
\frac{ \prod\limits_{\sigma\in I_0}(1-\tt^{\,\mm_{\,i(\sigma)}}
- \tt^{\,\mm_{\,j(\sigma)}}
+ \LL\; \tt^{\,\mm_{\,i(\sigma)}} \tt^{\,\mm_{\,j(\sigma)}})
}{\prod\limits_{i=1}^{s}
(1-\tt^{\,\mm_{\,i}})
(1-\LL\; \tt^{\,\mm_{\,i}})
}\; .
$$
\end{theorem}

\begin{proof}
One has the following equalities:
$$
\int_{\P \widehat S_D} \tt^{\,\ww}\ d \chi_g =
\int_{Y^D} \tt^{\,\ww}\ d \chi_g =
\sum\limits_{I\subset I_0}\quad \sum\limits_{\n\in {\cal N}^D_I}
[Y_{\n}]\cdot
\tt^{\,\vv(\n)}=
$$
$$
= \sum\limits_{I\subset I_0}\quad \sum\limits_{\n\in {\cal N}^D_I}
(\C^*)^{\# I} \prod\limits_{i=1}^{s}
[S^{n_i}\stackrel{\bullet}{E_i}]\cdot
\tt^{\ \sum\limits_{i=1}^r (n_i +
\sum\limits_{\sigma\in I :\ i(\sigma)=i}n'_{\sigma}   +
\sum\limits_{\sigma\in I :\ j(\sigma)=i}n''_{\sigma} )\, \mm_{\,i} } =
$$
$$
= \left( \sum\limits_{n_i\ge 0;\ i=1,\ldots,r}
\prod\limits_{i=1}^s
[S^{n_i}\stackrel{\bullet}{E_i}]\cdot
\tt^{\ \sum\limits_{i=1}^s n_i\mm_{\,i}}\right)
\cdot
\left(
\sum\limits_{I\subset I_0} (\C^*)^{\# I}
\sum
\limits_{n'_{\sigma}>0,\
n''_{\sigma}>0 :\ \sigma\in I}
\tt^{\,n'_{\sigma}\mm_{\,i(\sigma)}}\cdot
\tt^{\,n''_{\sigma}\mm_{\,j(\sigma)}}
\right).
$$
For the first factor one has
$$
\sum\limits_{n_i\ge 0;\ i=1,\ldots,r}
\prod\limits_{i=1}^s
[S^{n_i}\stackrel{\bullet}{E_i}]\cdot
\tt^{\ \sum\limits_{i=1}^s n_i\,\mm_{\,i}} =
 \prod\limits_{i=1}^s \left(
\sum\limits_{n=0}^{\infty}
[S^{n}\stackrel{\bullet}{E_i}]\cdot
\tt^{\,n\mm_{\,i}} \right)
=
\prod\limits_{i=1}^s
(1-\tt^{\,\mm_{\,i}}) ^{\,- [\stackrel{\bullet}{E_i}]} =
$$
(the last expression is in the sense of \cite{glm}; it looks like a motivic
version of the formula of A'Campo \cite{A'C})
$$
=
\prod\limits_{i=1}^s
(1-\tt^{\,\mm_{\,i}}) ^{- [\C \P^1]+ \# \{\sigma\in I_0:\ i(\sigma)=i\}
+ \# \{\sigma\in I_0 :\ j(\sigma)=i\}
} =
$$
$$
= \prod\limits_{i=1}^s
(1-\tt^{\,\mm_{\,i}}) ^{-(\LL +1)} \prod\limits_{\sigma\in I_0}
(1-\tt^{\,\mm_{\,i(\sigma)}})
(1-\tt^{\,\mm_{\,j(\sigma)}}) =
$$
\begin{equation}\label{eq12}
= \prod\limits_{i=1}^s
\frac{1}{
(1-\tt^{\,\mm_{\,i}}) (1-\LL \; \tt^{\,\mm_{\,i}})}
\prod\limits_{\sigma\in I_0}
(1-\tt^{\,\mm_{\,i(\sigma)}})
(1-\tt^{\,\mm_{\,j(\sigma)}})
\end{equation}
(the last equality follows from the fact that
$(1-t)^{-\LL} = 1/(1-\LL \; t)$).

\noindent
For the second factor one has
$$
\sum\limits_{I\subset I_0} (\C^*)^{\# I}
\sum\limits_{n'_{\sigma}>0,\
n''_{\sigma}>0 :\ \sigma\in I} \tt^{\,n'_{\sigma}\mm_{\,i(\sigma)}}\
\tt^{\,n''_{\sigma}\mm_{\,j(\sigma)}} =
$$
$$
=
\sum\limits_{I\subset I_0} (\C^*)^{\# I} \cdot
\prod\limits_{\sigma\in I}
\frac{\tt^{\,\mm_{\,i(\sigma)}}}{1-\tt^{\,\mm_{\,i(\sigma)}}}\cdot
\frac{\tt^{\,\mm_{\,j(\sigma)}}}{1-\tt^{\,\mm_{\,j(\sigma)}}}
 =
$$
\begin{equation}\label{eq13}
= \prod_{\sigma\in I_0} \left(
1+[\C^*]
\frac{\tt^{\,\mm_{\,i(\sigma)}}}{1-\tt^{\,\mm_{\,i(\sigma)}}}\cdot
\frac{\tt^{\,\mm_{\,j(\sigma)}}}{1-\tt^{\,\mm_{\,j(\sigma)}}}
\right)\ .
\end{equation}
Combining (\ref{eq12}) and (\ref{eq13}) one gets
$$
\int\limits_{\P \widehat S_D} \tt^{\,\vv}\ d \chi_g =
\prod\limits_{i=1}^s
\frac{1}{
(1-\tt^{\,\mm_{\,i}}) (1-\LL \; \tt^{\,\mm_{\,i}})}
\prod\limits_{\sigma\in I_0}
\left[ (1-\tt^{\,\mm_{\,i(\sigma)}})
(1-\tt^{\,\mm_{\,j(\sigma)}})
+ [\C^*]\,\tt^{\,\mm_{\,i(\sigma)}}
\tt^{\,\mm_{\,j(\sigma)}}
\right] =
$$
$$
=
\frac{ \prod\limits_{\sigma\in I_0}(1-\tt^{\,\mm_{\,i(\sigma)}}
- \tt^{\,\mm_{\,j(\sigma)}}
+ \LL\; \tt^{\,\mm_{\,i(\sigma)}} \tt^{\,\mm_{\,j(\sigma)}})
}{\prod\limits_{i=1}^{s}
(1-\tt^{\,\mm_{\,i}})
(1-\LL\; \tt^{\,\mm_{\,i}})
}\; .
$$
\end{proof}

\end{document}